\documentclass[12pt]{amsart}
\usepackage{amssymb,mathtools}
\usepackage[T1]{fontenc}
\usepackage{lmodern}
\usepackage[english]{babel}
\usepackage{cite}
\usepackage{xurl}
\usepackage{needspace}
\usepackage{tikz}
\usepackage[hidelinks]{hyperref}

\newtheorem{proposition}{Proposition}[section]
\newtheorem{theorem}[proposition]{Theorem}
\newtheorem{maintheorem}[proposition]{Main Theorem}
\newtheorem{lemma}[proposition]{Lemma}
\newtheorem{corollary}[proposition]{Corollary}
\theoremstyle{definition}
\newtheorem{definition}[proposition]{Definition}
\theoremstyle{remark}
\newtheorem{remark}[proposition]{Remark}
\newcommand{\inv}{^{-1}}

\begin{document}
\title[Parabolic double cosets]{The double coset problem for parabolic subgroups of braid groups}
\author{Arkadius Kalka}
\author{Mina Teicher}
\author{Boaz Tsaban}
\address{Department of Mathematics, Bar-Ilan University, Ramat Gan, Israel}
\email{Arkadius.Kalka@rub.de}
\email{teicher@math.biu.ac.il}
\email{tsaban@math.biu.ac.il}

\begin{abstract}
We solve the double coset problem for arbitrary parabolic subgroups of the $n$-strand braid group
$B_n$, by reducing it to the simultaneous conjugacy problem in the larger group $B_{n+1}$.
The resulting decision algorithm is uniform in the braid index
and the supplied subgroup generators.
We also formulate the reduction for general
group embeddings, and apply it to additional commuting constraints,
several other spherical Artin types, and finite direct products.
\end{abstract}
\subjclass[2020]{20F36; 20F10}
\keywords{braid groups, parabolic subgroups, double coset problem, simultaneous conjugacy, double centralizers}
\maketitle
\raggedbottom

\section{The double coset problem}\label{sec:introduction}

Decision problems play an important role in combinatorial group
theory~[\citen{LyndonSchupp1977}; \citen{Miller1992}].
Elements of a finitely generated group are represented by words in a
specified finite generating set and its inverses. The \emph{word problem}
asks whether a given word represents the identity element. Given two
finite lists $(u_1,\ldots,u_k)$ and $(v_1,\ldots,v_k)$ of elements of the
group, the
\emph{simultaneous conjugacy problem} asks whether there is a group
element $x$ satisfying $x\inv u_i x=v_i$ for all $i=1,\dotsc,k$.

For each positive integer $n$, let $B_n$ denote the braid group on $n$
strands. Geometrically, braids are disjoint strands running monotonically
between two fixed rows of endpoints, considered up to deformation fixing
the endpoints; multiplication is stacking.
For $i=1,\dotsc,n-1$, the standard Artin generator $\sigma_i$
makes one positive crossing in which strand $i+1$ passes over strand $i$,
leaving the other strands straight. Its inverse reverses that
crossing. The pure braid group $P_n$ consists of braids in which
each strand returns to its original position.
Figure~\ref{fig:basic-braids} illustrates these basic braid operations.

\begin{figure}[htbp]
\centering
\begin{tikzpicture}[
  blue strand/.style={draw=blue!65!black,line width=1pt},
  orange strand/.style={draw=orange!80!black,line width=1pt},
  gray strand/.style={draw=black!35,line width=.7pt},
  endpoint label/.style={font=\scriptsize,inner sep=1pt},
  braid title/.style={font=\small,align=center,text width=3.4cm}
]
\node[braid title] at (1.3,2.93) {Generator $\sigma_1$};
\node[braid title] at (5,2.93) {Inverse braid $\sigma_1^{-1}$};
\node[braid title] at (8.7,2.93) {Pure braid $\sigma_1^2$};
\draw[blue strand] (0.7,2.2)
  .. controls (0.7,1.4667) and (1.3,0.7333) .. (1.3,0);
\draw[orange strand] (1.3,2.2)
  .. controls (1.3,1.4667) and (0.7,0.7333) .. (0.7,0);
\draw[gray strand] (1.9,2.2)
  .. controls (1.9,1.4667) and (1.9,0.7333) .. (1.9,0);
\draw[orange strand,preaction={draw=white,line width=3.1pt}] (1.1779,1.562) .. controls (1.0741,1.254) and (0.9259,0.946) .. (0.8221,0.638);
\node[endpoint label] at (0.7,2.39) {$1$};
\node[endpoint label] at (0.7,-0.19) {$2$};
\node[endpoint label] at (1.3,2.39) {$2$};
\node[endpoint label] at (1.3,-0.19) {$1$};
\node[endpoint label] at (1.9,2.39) {$3$};
\node[endpoint label] at (1.9,-0.19) {$3$};
\draw[blue strand] (4.4,2.2)
  .. controls (4.4,1.4667) and (5,0.7333) .. (5,0);
\draw[orange strand] (5,2.2)
  .. controls (5,1.4667) and (4.4,0.7333) .. (4.4,0);
\draw[gray strand] (5.6,2.2)
  .. controls (5.6,1.4667) and (5.6,0.7333) .. (5.6,0);
\draw[blue strand,preaction={draw=white,line width=3.1pt}] (4.5221,1.562) .. controls (4.6259,1.254) and (4.7741,0.946) .. (4.8779,0.638);
\node[endpoint label] at (4.4,2.39) {$1$};
\node[endpoint label] at (4.4,-0.19) {$2$};
\node[endpoint label] at (5,2.39) {$2$};
\node[endpoint label] at (5,-0.19) {$1$};
\node[endpoint label] at (5.6,2.39) {$3$};
\node[endpoint label] at (5.6,-0.19) {$3$};
\draw[blue strand] (8.1,2.2)
  .. controls (8.1,1.8333) and (8.7,1.4667) .. (8.7,1.1)
  .. controls (8.7,0.7333) and (8.1,0.3667) .. (8.1,0);
\draw[orange strand] (8.7,2.2)
  .. controls (8.7,1.8333) and (8.1,1.4667) .. (8.1,1.1)
  .. controls (8.1,0.7333) and (8.7,0.3667) .. (8.7,0);
\draw[gray strand] (9.3,2.2)
  .. controls (9.3,1.8333) and (9.3,1.4667) .. (9.3,1.1)
  .. controls (9.3,0.7333) and (9.3,0.3667) .. (9.3,0);
\draw[orange strand,preaction={draw=white,line width=3.1pt}] (8.5779,1.881) .. controls (8.4741,1.727) and (8.3259,1.573) .. (8.2221,1.419);
\draw[blue strand,preaction={draw=white,line width=3.1pt}] (8.5779,0.781) .. controls (8.4741,0.627) and (8.3259,0.473) .. (8.2221,0.319);
\node[endpoint label] at (8.1,2.39) {$1$};
\node[endpoint label] at (8.1,-0.19) {$1$};
\node[endpoint label] at (8.7,2.39) {$2$};
\node[endpoint label] at (8.7,-0.19) {$2$};
\node[endpoint label] at (9.3,2.39) {$3$};
\node[endpoint label] at (9.3,-0.19) {$3$};
\end{tikzpicture}
\caption{Basic braids in $B_3$}
\label{fig:basic-braids}
\end{figure}

The word and conjugacy problems in braid groups were first solved by
Artin~[\citen{Ar25}, Section~7; \citen{Ar47}, Theorem~14] and Garside~\cite[Theorem~10]{Ga69}, respectively.
Lee and Lee~\cite[Theorems~2 and~3]{LL02} solved the simultaneous conjugacy problem in braid
groups. Kalka, Tsaban and Vinokur~\cite[Sections~6--7]{KTV14} refined their solution.

\begin{definition}[Double Coset Problem]\label{def:double-coset}
Given subgroups $A,B$ of a group $G$ and group elements $g,h\in G$,
decide whether $h\in AgB$, that is, whether
$h=agb$ for some group elements $a\in A$ and $b\in B$.
The associated search problem is to find such factors, provided that they exist.
\end{definition}

Subgroups are specified by finite lists of generators, written as words
in the ambient generators and their inverses.
The \emph{subgroup membership problem} is to decide, given
a subgroup $H$ and a group element $w$, whether $w\in H$.
This is the special case $w\in H\cdot1\cdot\{1\}$ of the double coset problem.

Mikhailova~[\citen{Mi58}, Theorem~1; \citen{Mi66}, Theorem~1] constructed a finitely generated subgroup of
$F_2\times F_2$ with undecidable membership, where $F_2$ is the free group
of rank two. Using an embedding of this direct product in the braid group
$B_n$, Makanina~\cite[Theorems~1 and~2]{Ma81} established subgroup-membership undecidability
for braid indices $n\ge5$.
We will show that for parabolic subgroups, this problem is solvable.

\begin{definition}\label{def:parabolic}
A subgroup of the braid group $B_n$ is \emph{standard parabolic} if it is
generated by a subset of the standard Artin generators
$\{\sigma_1,\ldots,\sigma_{n-1}\}$.
A subgroup is \emph{parabolic} if it is conjugate to a standard parabolic
subgroup. Geometrically, a standard parabolic subgroup allows independent
braiding within disjoint consecutive blocks of strands.
\end{definition}

\begin{figure}[htbp]
\centering
\begin{tikzpicture}[
  blue strand/.style={draw=blue!65!black,line width=1pt},
  orange strand/.style={draw=orange!80!black,line width=1pt},
  blue block/.style={draw=blue!35!black,fill=blue!5,line width=.4pt,rounded corners=2pt},
  orange block/.style={draw=orange!55!black,fill=orange!7,line width=.4pt,rounded corners=2pt},
  endpoint label/.style={font=\scriptsize,inner sep=1pt},
  braid title/.style={font=\small,align=center,text width=4.8cm},
  braid label/.style={font=\footnotesize,align=center,text width=1.45cm},
  arrow label/.style={font=\footnotesize,align=center}
]
\node[braid title] at (1.5,2.95) {Standard\\parabolic subgroup};
\node[braid title] at (7.8,2.95) {Parabolic subgroup};
\path[blue block] (.36,-.05) rectangle (1.44,2.25);
\path[orange block] (1.56,-.05) rectangle (2.64,2.25);
\path[blue block] (6.66,.55) rectangle (7.74,1.65);
\path[orange block] (7.86,.55) rectangle (8.94,1.65);
\draw[->,line width=.65pt] (3.5,1.1) -- (5.65,1.1);
\node[arrow label] at (4.575,1.6) {$\beta\mapsto\sigma_2^{-1}\beta\sigma_2$};
\node[braid label] at (9.75,1.925) {$\sigma_2^{-1}$};
\node[braid label] at (9.75,.275) {$\sigma_2$};
\draw[blue strand] (0.6,2.2)
  .. controls (0.6,1.8333) and (1.2,1.4667) .. (1.2,1.1)
  .. controls (1.2,0.7333) and (1.2,0.3667) .. (1.2,0);
\draw[blue strand] (1.2,2.2)
  .. controls (1.2,1.8333) and (0.6,1.4667) .. (0.6,1.1)
  .. controls (0.6,0.7333) and (0.6,0.3667) .. (0.6,0);
\draw[orange strand] (1.8,2.2)
  .. controls (1.8,1.8333) and (1.8,1.4667) .. (1.8,1.1)
  .. controls (1.8,0.7333) and (2.4,0.3667) .. (2.4,0);
\draw[orange strand] (2.4,2.2)
  .. controls (2.4,1.8333) and (2.4,1.4667) .. (2.4,1.1)
  .. controls (2.4,0.7333) and (1.8,0.3667) .. (1.8,0);
\draw[blue strand,preaction={draw=white,line width=3.1pt}] (0.7221,1.881) .. controls (0.8259,1.727) and (0.9741,1.573) .. (1.0779,1.419);
\draw[orange strand,preaction={draw=white,line width=3.1pt}] (1.9221,0.781) .. controls (2.0259,0.627) and (2.1741,0.473) .. (2.2779,0.319);
\node[endpoint label] at (0.6,2.39) {$1$};
\node[endpoint label] at (0.6,-0.19) {$2$};
\node[endpoint label] at (1.2,2.39) {$2$};
\node[endpoint label] at (1.2,-0.19) {$1$};
\node[endpoint label] at (1.8,2.39) {$3$};
\node[endpoint label] at (1.8,-0.19) {$4$};
\node[endpoint label] at (2.4,2.39) {$4$};
\node[endpoint label] at (2.4,-0.19) {$3$};
\draw[blue strand] (6.9,2.2)
  .. controls (6.9,2.0167) and (6.9,1.8333) .. (6.9,1.65)
  .. controls (6.9,1.4667) and (7.5,1.2833) .. (7.5,1.1)
  .. controls (7.5,0.9167) and (7.5,0.7333) .. (7.5,0.55)
  .. controls (7.5,0.3667) and (8.1,0.1833) .. (8.1,0);
\draw[orange strand] (7.5,2.2)
  .. controls (7.5,2.0167) and (8.1,1.8333) .. (8.1,1.65)
  .. controls (8.1,1.4667) and (8.1,1.2833) .. (8.1,1.1)
  .. controls (8.1,0.9167) and (8.7,0.7333) .. (8.7,0.55)
  .. controls (8.7,0.3667) and (8.7,0.1833) .. (8.7,0);
\draw[blue strand] (8.1,2.2)
  .. controls (8.1,2.0167) and (7.5,1.8333) .. (7.5,1.65)
  .. controls (7.5,1.4667) and (6.9,1.2833) .. (6.9,1.1)
  .. controls (6.9,0.9167) and (6.9,0.7333) .. (6.9,0.55)
  .. controls (6.9,0.3667) and (6.9,0.1833) .. (6.9,0);
\draw[orange strand] (8.7,2.2)
  .. controls (8.7,2.0167) and (8.7,1.8333) .. (8.7,1.65)
  .. controls (8.7,1.4667) and (8.7,1.2833) .. (8.7,1.1)
  .. controls (8.7,0.9167) and (8.1,0.7333) .. (8.1,0.55)
  .. controls (8.1,0.3667) and (7.5,0.1833) .. (7.5,0);
\draw[orange strand,preaction={draw=white,line width=3.1pt}] (7.6221,2.0405) .. controls (7.7259,1.9635) and (7.8741,1.8865) .. (7.9779,1.8095);
\draw[blue strand,preaction={draw=white,line width=3.1pt}] (7.0221,1.4905) .. controls (7.1259,1.4135) and (7.2741,1.3365) .. (7.3779,1.2595);
\draw[orange strand,preaction={draw=white,line width=3.1pt}] (8.2221,0.9405) .. controls (8.3259,0.8635) and (8.4741,0.7865) .. (8.5779,0.7095);
\draw[orange strand,preaction={draw=white,line width=3.1pt}] (7.9779,0.3905) .. controls (7.8741,0.3135) and (7.7259,0.2365) .. (7.6221,0.1595);
\node[endpoint label] at (6.9,2.39) {$1$};
\node[endpoint label] at (6.9,-0.19) {$3$};
\node[endpoint label] at (7.5,2.39) {$2$};
\node[endpoint label] at (7.5,-0.19) {$4$};
\node[endpoint label] at (8.1,2.39) {$3$};
\node[endpoint label] at (8.1,-0.19) {$1$};
\node[endpoint label] at (8.7,2.39) {$4$};
\node[endpoint label] at (8.7,-0.19) {$2$};
\end{tikzpicture}
\caption{Elements of parabolic subgroups}
\label{fig:parabolic-braids}
\end{figure}

Figure~\ref{fig:parabolic-braids} illustrates a standard parabolic
subgroup and its conjugate.
For Coxeter groups, standard parabolic double cosets have unique representatives
of minimal length~\cite[Chapter~IV, \S1, Exercise~3]{Bo68}. Digne, Godelle and
Michel~\cite[Proposition~0.3]{DGM24} obtain a related partial result for
Artin groups using parabolic retractions. We treat arbitrary pairs of
parabolic subgroups of braid groups, without a connectedness assumption
on their Coxeter diagrams or a restriction on the double coset
representative.

\begin{maintheorem}\label{thm:main}
There is an algorithm that, given a positive integer $n$, finite
generating sets for parabolic subgroups $A,B\le B_n$, and braids
$g,h\in B_n$, decides whether $h\in AgB$.
For each positive instance, it produces words $a\in A$ and $b\in B$, in the corresponding generators,
such that $h=agb$.
\end{maintheorem}

Parabolicity is a promise on the input subgroups. The algorithm does not
require their standard parabolic forms or conjugating braids.
For each braid index $n\ge2$, the decision problem reduces
to a single simultaneous conjugacy instance in the braid group $B_{n+1}$.

Section~\ref{sec:criterion} provides a general reduction of the double coset
problem to simultaneous conjugacy in a larger group.
Section~\ref{sec:braids} applies this reduction to braid groups by adding
a trivial strand, proves the main theorem, and estimates the size of
the resulting simultaneous conjugacy instance.
Section~\ref{sec:discussion} generalizes the main theorem to some other spherical Artin
groups, and to finite direct products.

\section{A reduction through a larger group}\label{sec:criterion}

For a subset $S$ of a group $G$, let
\[
C_G(S)\coloneq\{x\in G:xs=sx\text{ for all }s\in S\}
\]
be the centralizer of the set $S$ in the group $G$,
and $Z(G)\coloneq C_G(G)$ be the center of the group.
If $S$ generates $H$, then $C_G(S)=C_G(H)$. The subgroup
$C_G(C_G(S))$ is the \emph{double centralizer} of $S$.

\begin{proposition}\label{prop:ambient-reduction}
Let $H\le G$ be finitely generated groups, with the inclusion specified
by words for the generators of $H$ in those of $G$.
Let $A,B\le H$ be finitely generated subgroups such that
$C_G(C_G(A))=Z(G)A$ and $C_G(C_G(B))=Z(G)B$, and
finite generating sets for $C_G(A)$ and $C_G(B)$ can be computed.
If $Z(G)\cap H=\{1\}$, then the double coset problem for the subgroups
$A,B$ of $H$ effectively reduces to a simultaneous conjugacy
instance in $G$. The reduction is uniform whenever the centralizer
computations are uniform.
\end{proposition}

\begin{proof}
Compute
\begin{align*}
C_G(A)&=\langle c_1,\ldots,c_r\rangle,\\
C_G(B)&=\langle d_1,\ldots,d_s\rangle.
\end{align*}
For $g,h\in H$, consider the following system, with unknown $x\in G$:
\begin{equation}\label{eq:simultaneous}
\begin{cases}
x\inv c_i x=c_i & 1\le i\le r,\\
x\inv hd_jh\inv x=gd_jg\inv & 1\le j\le s.
\end{cases}
\end{equation}
All entries are computable words in the generators of $G$.
If $h=agb$ with $a\in A$ and $b\in B$, then $x=a$ solves the system:
The element $a$ commutes with the elements $c_i$, and
\[
a\inv hd_jh\inv a = a\inv (agb)d_j(agb)\inv a =gbd_jb\inv g\inv=gd_jg\inv.
\]

Conversely, let $x$ be a solution of the system, and set
$y\coloneq g\inv x\inv h$. The first family of equations implies
$x\in C_G(C_G(A))=Z(G)A$. The second gives
\[
yd_jy\inv=g\inv x\inv hd_jh\inv xg=d_j,
\]
so $y\in C_G(C_G(B))=Z(G)B$.
Write $x=z_Aa$ and $y=z_Bb$, where $z_A,z_B\in Z(G)$, $a\in A$, and
$b\in B$. Since $h=xgy$, we have
\[
h=(z_Az_B)agb,\qquad
z_Az_B=h(agb)\inv\in Z(G)\cap H=\{1\}.
\]
Thus $h=agb$.
\end{proof}

Without the condition $Z(G)\cap H=\{1\}$, the same system is solvable
exactly when $h\in Z(G)AgB$. Indeed, the converse calculation still
applies, and if $h=zagb$ with $z\in Z(G)$, then $x=a$ solves the system
because $z$ disappears under conjugation. The intersection condition
therefore removes the possible discrepancy between
$Z(G)AgB\cap H$ and $AgB$.

\begin{remark}\label{rem:witnesses}
If the word problem in $H$ is decidable, a positive decision answer
suffices to recover factors in the supplied subgroup generators.
Enumerate pairs of words in the generators of $A$ and $B$ and their
inverses in increasing total length, including the empty words, and test
$h=agb$ for each pair. Each stage is finite, and a witnessing pair is
eventually reached. This gives a uniform search algorithm whenever the
word-problem and decision algorithms are uniform.
\end{remark}

\section{Proof of the Main Theorem}\label{sec:braids}

Centralizers of parabolic subgroups of braid groups were studied by
Fenn, Rolfsen and Zhu, and Paris~[\citen{FRZ96}; \citen{Pa97}].
The following double centralizer theorem was established by Garber,
Kalka, Liberman and Teicher~\cite[Theorem~3.16]{GKLT13}.
It also follows from the spherical-type double centralizer theorem of
Ajbal and Godelle~\cite[Theorem~0.1 and Proposition~2.1]{AG19}: the exceptional cases do not occur in
Coxeter type $A$ of braid groups, and their theorem allows reducible
standard parabolic subgroups. Conjugation gives the assertion for nonstandard parabolic
subgroups. The cases of the trivial subgroup and the whole group are
immediate.

The braid $\Delta_n\in B_n$ is the Garside half twist, which crosses each
pair of strands once positively. Its square $\Delta_n^2$ is the full
twist, in which each pair of strands winds once around one another.

\begin{theorem}\label{thm:double-centralizer}
For each $n\ge3$ and every parabolic subgroup $H\le B_n$, we have
\[
C_{B_n}(C_{B_n}(H))=Z(B_n)H
=\langle\Delta_n^2\rangle H.
\]
\end{theorem}

Gurzo's algorithm~\cite[Theorem]{Gu85} computes a finite generating set for the
centralizer of any finite set of braids.
Lee and Lee's algorithm~\cite[Theorems~2 and~3]{LL02} decides simultaneous conjugacy in
braid groups. We use these algorithms in their uniform forms: the strand
number and the relevant finite lists of braid words are part of the input.
Similarly, the braid word-problem algorithm is used in its uniform form.
When the input generating set is empty, its centralizer is the whole
braid group, whose standard generators are already available.

For each braid index $n\ge2$, the standard embedding $B_n \hookrightarrow B_{n+1}$ identifies the braid group $B_n$ with the
subgroup
\[
B_n\coloneq\langle\sigma_1,\dotsc,\sigma_{n-1}\rangle\le B_{n+1}
\]
by adding a trivial last strand.
This embedding is effective: it leaves each word in the original
standard generators unchanged.

\begin{lemma}\label{lem:added-strand}
For each braid index $n\ge2$, the standard embedding $B_n \hookrightarrow B_{n+1}$ satisfies
\[
Z(B_{n+1})\cap B_n=\{1\}.
\]
Moreover, every parabolic subgroup of the braid group $B_n$ is parabolic
in the larger braid group $B_{n+1}$ under this embedding.
\end{lemma}

\begin{proof}
The center of the braid group $B_{n+1}$ is generated by the full twist
$\Delta_{n+1}^2$~\cite[Theorem~7(ii)]{Ga69}.
Suppose that a group element $z\coloneq\Delta_{n+1}^{2q}$ belongs to the subgroup
$B_n$, where the exponent $q$ is an integer.
Since the braid $z$ is pure, we may apply the strand-deletion homomorphism
$P_{n+1}\to P_2$ that retains only the first and last strands. Purity
ensures that deleting a fixed set of strands respects stacking.

Every braid in the subgroup $B_n\cap P_{n+1}$ maps to the identity: its
added last strand is trivial and has no crossings with the original
strands. The full twist maps to the two-strand full twist $\sigma_1^2$ (see Figure~\ref{fig:strand-deletion}).
Thus, we have $\sigma_1^{2q}=1$.
Since the two-strand braid group $B_2$ is infinite
cyclic with generator $\sigma_1$, the exponent $q$ is zero.
This proves the intersection assertion.
The second assertion is immediate.
\end{proof}

\begin{figure}[htbp]
\centering
\begin{tikzpicture}[
  first strand/.style={draw=blue!65!black,line width=1pt},
  last strand/.style={draw=orange!80!black,line width=1pt},
  deleted strand/.style={draw=black!35,line width=.7pt},
  endpoint label/.style={font=\scriptsize,inner sep=1pt},
  braid title/.style={font=\small,align=center,text width=4.7cm},
  deletion arrow/.style={->,line width=.65pt},
  deletion label/.style={font=\footnotesize,align=center}
]
\node[braid title] at (1.7,3.35) {(a) Embedded braid $\sigma_1^2$};
\node[braid title] at (8.15,3.35) {Identity braid};
\draw[deletion arrow] (3.6,1.3) -- (6.55,1.3);
\node[deletion label] at (5.075,1.75) {delete strands\\2 and 3};
\draw[first strand] (0.95,2.6)
  .. controls (0.95,2.1667) and (1.45,1.7333) .. (1.45,1.3)
  .. controls (1.45,0.8667) and (0.95,0.4333) .. (0.95,0);
\draw[deleted strand] (1.45,2.6)
  .. controls (1.45,2.1667) and (0.95,1.7333) .. (0.95,1.3)
  .. controls (0.95,0.8667) and (1.45,0.4333) .. (1.45,0);
\draw[deleted strand] (1.95,2.6)
  .. controls (1.95,2.1667) and (1.95,1.7333) .. (1.95,1.3)
  .. controls (1.95,0.8667) and (1.95,0.4333) .. (1.95,0);
\draw[last strand] (2.45,2.6)
  .. controls (2.45,2.1667) and (2.45,1.7333) .. (2.45,1.3)
  .. controls (2.45,0.8667) and (2.45,0.4333) .. (2.45,0);
\draw[deleted strand,preaction={draw=white,line width=3.1pt}] (1.3482,2.223) .. controls (1.2618,2.041) and (1.1382,1.859) .. (1.0518,1.677);
\draw[first strand,preaction={draw=white,line width=3.1pt}] (1.3482,0.923) .. controls (1.2618,0.741) and (1.1382,0.559) .. (1.0518,0.377);
\node[endpoint label] at (0.95,2.79) {$1$};
\node[endpoint label] at (0.95,-0.19) {$1$};
\node[endpoint label] at (1.45,2.79) {$2$};
\node[endpoint label] at (1.45,-0.19) {$2$};
\node[endpoint label] at (1.95,2.79) {$3$};
\node[endpoint label] at (1.95,-0.19) {$3$};
\node[endpoint label] at (2.45,2.79) {$4$};
\node[endpoint label] at (2.45,-0.19) {$4$};
\draw[first strand] (7.65,2.6)
  .. controls (7.65,1.7333) and (7.65,0.8667) .. (7.65,0);
\draw[last strand] (8.65,2.6)
  .. controls (8.65,1.7333) and (8.65,0.8667) .. (8.65,0);
\node[endpoint label] at (7.65,2.79) {$1$};
\node[endpoint label] at (7.65,-0.19) {$1$};
\node[endpoint label] at (8.65,2.79) {$4$};
\node[endpoint label] at (8.65,-0.19) {$4$};
\begin{scope}[shift={(0,-4.05)}]
\node[braid title] at (1.7,3.35) {(b) Full twist $\Delta_4^2$};
\node[braid title] at (8.15,3.35) {Two-strand full twist $\sigma_1^2$};
\draw[deletion arrow] (3.6,1.3) -- (6.55,1.3);
\node[deletion label] at (5.075,1.75) {delete strands\\2 and 3};
\draw[first strand] (0.95,2.6)
  .. controls (0.95,2.5278) and (1.45,2.4556) .. (1.45,2.3833)
  .. controls (1.45,2.3111) and (1.95,2.2389) .. (1.95,2.1667)
  .. controls (1.95,2.0944) and (2.45,2.0222) .. (2.45,1.95)
  .. controls (2.45,1.8778) and (2.45,1.8056) .. (2.45,1.7333)
  .. controls (2.45,1.6611) and (2.45,1.5889) .. (2.45,1.5167)
  .. controls (2.45,1.4444) and (1.95,1.3722) .. (1.95,1.3)
  .. controls (1.95,1.2278) and (1.95,1.1556) .. (1.95,1.0833)
  .. controls (1.95,1.0111) and (1.45,0.9389) .. (1.45,0.8667)
  .. controls (1.45,0.7944) and (1.45,0.7222) .. (1.45,0.65)
  .. controls (1.45,0.5778) and (0.95,0.5056) .. (0.95,0.4333)
  .. controls (0.95,0.3611) and (0.95,0.2889) .. (0.95,0.2167)
  .. controls (0.95,0.1444) and (0.95,0.0722) .. (0.95,0);
\draw[deleted strand] (1.45,2.6)
  .. controls (1.45,2.5278) and (0.95,2.4556) .. (0.95,2.3833)
  .. controls (0.95,2.3111) and (0.95,2.2389) .. (0.95,2.1667)
  .. controls (0.95,2.0944) and (0.95,2.0222) .. (0.95,1.95)
  .. controls (0.95,1.8778) and (1.45,1.8056) .. (1.45,1.7333)
  .. controls (1.45,1.6611) and (1.95,1.5889) .. (1.95,1.5167)
  .. controls (1.95,1.4444) and (2.45,1.3722) .. (2.45,1.3)
  .. controls (2.45,1.2278) and (2.45,1.1556) .. (2.45,1.0833)
  .. controls (2.45,1.0111) and (2.45,0.9389) .. (2.45,0.8667)
  .. controls (2.45,0.7944) and (1.95,0.7222) .. (1.95,0.65)
  .. controls (1.95,0.5778) and (1.95,0.5056) .. (1.95,0.4333)
  .. controls (1.95,0.3611) and (1.45,0.2889) .. (1.45,0.2167)
  .. controls (1.45,0.1444) and (1.45,0.0722) .. (1.45,0);
\draw[deleted strand] (1.95,2.6)
  .. controls (1.95,2.5278) and (1.95,2.4556) .. (1.95,2.3833)
  .. controls (1.95,2.3111) and (1.45,2.2389) .. (1.45,2.1667)
  .. controls (1.45,2.0944) and (1.45,2.0222) .. (1.45,1.95)
  .. controls (1.45,1.8778) and (0.95,1.8056) .. (0.95,1.7333)
  .. controls (0.95,1.6611) and (0.95,1.5889) .. (0.95,1.5167)
  .. controls (0.95,1.4444) and (0.95,1.3722) .. (0.95,1.3)
  .. controls (0.95,1.2278) and (1.45,1.1556) .. (1.45,1.0833)
  .. controls (1.45,1.0111) and (1.95,0.9389) .. (1.95,0.8667)
  .. controls (1.95,0.7944) and (2.45,0.7222) .. (2.45,0.65)
  .. controls (2.45,0.5778) and (2.45,0.5056) .. (2.45,0.4333)
  .. controls (2.45,0.3611) and (2.45,0.2889) .. (2.45,0.2167)
  .. controls (2.45,0.1444) and (1.95,0.0722) .. (1.95,0);
\draw[last strand] (2.45,2.6)
  .. controls (2.45,2.5278) and (2.45,2.4556) .. (2.45,2.3833)
  .. controls (2.45,2.3111) and (2.45,2.2389) .. (2.45,2.1667)
  .. controls (2.45,2.0944) and (1.95,2.0222) .. (1.95,1.95)
  .. controls (1.95,1.8778) and (1.95,1.8056) .. (1.95,1.7333)
  .. controls (1.95,1.6611) and (1.45,1.5889) .. (1.45,1.5167)
  .. controls (1.45,1.4444) and (1.45,1.3722) .. (1.45,1.3)
  .. controls (1.45,1.2278) and (0.95,1.1556) .. (0.95,1.0833)
  .. controls (0.95,1.0111) and (0.95,0.9389) .. (0.95,0.8667)
  .. controls (0.95,0.7944) and (0.95,0.7222) .. (0.95,0.65)
  .. controls (0.95,0.5778) and (1.45,0.5056) .. (1.45,0.4333)
  .. controls (1.45,0.3611) and (1.95,0.2889) .. (1.95,0.2167)
  .. controls (1.95,0.1444) and (2.45,0.0722) .. (2.45,0);
\draw[deleted strand,preaction={draw=white,line width=3.1pt}] (1.3482,2.5372) .. controls (1.2618,2.5068) and (1.1382,2.4765) .. (1.0518,2.4462);
\draw[deleted strand,preaction={draw=white,line width=3.1pt}] (1.8482,2.3205) .. controls (1.7618,2.2902) and (1.6382,2.2598) .. (1.5518,2.2295);
\draw[last strand,preaction={draw=white,line width=3.1pt}] (2.3482,2.1038) .. controls (2.2618,2.0735) and (2.1382,2.0432) .. (2.0518,2.0128);
\draw[deleted strand,preaction={draw=white,line width=3.1pt}] (1.3482,1.8872) .. controls (1.2618,1.8568) and (1.1382,1.8265) .. (1.0518,1.7962);
\draw[last strand,preaction={draw=white,line width=3.1pt}] (1.8482,1.6705) .. controls (1.7618,1.6402) and (1.6382,1.6098) .. (1.5518,1.5795);
\draw[first strand,preaction={draw=white,line width=3.1pt}] (2.3482,1.4538) .. controls (2.2618,1.4235) and (2.1382,1.3932) .. (2.0518,1.3628);
\draw[last strand,preaction={draw=white,line width=3.1pt}] (1.3482,1.2372) .. controls (1.2618,1.2068) and (1.1382,1.1765) .. (1.0518,1.1462);
\draw[first strand,preaction={draw=white,line width=3.1pt}] (1.8482,1.0205) .. controls (1.7618,0.9902) and (1.6382,0.9598) .. (1.5518,0.9295);
\draw[deleted strand,preaction={draw=white,line width=3.1pt}] (2.3482,0.8038) .. controls (2.2618,0.7735) and (2.1382,0.7432) .. (2.0518,0.7128);
\draw[first strand,preaction={draw=white,line width=3.1pt}] (1.3482,0.5872) .. controls (1.2618,0.5568) and (1.1382,0.5265) .. (1.0518,0.4962);
\draw[deleted strand,preaction={draw=white,line width=3.1pt}] (1.8482,0.3705) .. controls (1.7618,0.3402) and (1.6382,0.3098) .. (1.5518,0.2795);
\draw[deleted strand,preaction={draw=white,line width=3.1pt}] (2.3482,0.1538) .. controls (2.2618,0.1235) and (2.1382,0.0932) .. (2.0518,0.0628);
\node[endpoint label] at (0.95,2.79) {$1$};
\node[endpoint label] at (0.95,-0.19) {$1$};
\node[endpoint label] at (1.45,2.79) {$2$};
\node[endpoint label] at (1.45,-0.19) {$2$};
\node[endpoint label] at (1.95,2.79) {$3$};
\node[endpoint label] at (1.95,-0.19) {$3$};
\node[endpoint label] at (2.45,2.79) {$4$};
\node[endpoint label] at (2.45,-0.19) {$4$};
\draw[first strand] (7.65,2.6)
  .. controls (7.65,2.1667) and (8.65,1.7333) .. (8.65,1.3)
  .. controls (8.65,0.8667) and (7.65,0.4333) .. (7.65,0);
\draw[last strand] (8.65,2.6)
  .. controls (8.65,2.1667) and (7.65,1.7333) .. (7.65,1.3)
  .. controls (7.65,0.8667) and (8.65,0.4333) .. (8.65,0);
\draw[last strand,preaction={draw=white,line width=3.1pt}] (8.4465,2.223) .. controls (8.2735,2.041) and (8.0265,1.859) .. (7.8535,1.677);
\draw[first strand,preaction={draw=white,line width=3.1pt}] (8.4465,0.923) .. controls (8.2735,0.741) and (8.0265,0.559) .. (7.8535,0.377);
\node[endpoint label] at (7.65,2.79) {$1$};
\node[endpoint label] at (7.65,-0.19) {$1$};
\node[endpoint label] at (8.65,2.79) {$4$};
\node[endpoint label] at (8.65,-0.19) {$4$};
\end{scope}
\end{tikzpicture}
\caption{Strand deletion}
\label{fig:strand-deletion}
\end{figure}

\begin{proof}[Proof of Main Theorem~\ref{thm:main}]
The assertion is trivial for braid index $n=1$.
For braid indices $n\ge2$, use the effective embedding
$B_n\le B_{n+1}$.
By Lemma~\ref{lem:added-strand}, the input subgroups $A$ and $B$ are
parabolic in the larger group $B_{n+1}$, and $Z(B_{n+1})\cap B_n=\{1\}$.
By Theorem~\ref{thm:double-centralizer}, we have both double centralizer
identities required in Proposition~\ref{prop:ambient-reduction}.
Apply Gurzo's algorithm to the supplied generators of the subgroups $A$
and $B$, viewed as braids in the group $B_{n+1}$, to compute generating
sets for their centralizers $C_{B_{n+1}}(A)$ and $C_{B_{n+1}}(B)$.
Use Proposition~\ref{prop:ambient-reduction} to define a
simultaneous conjugacy instance in the group $B_{n+1}$.
Use the Lee--Lee algorithm to decide this instance.

The braid word problem is decidable uniformly in $n$, so
Remark~\ref{rem:witnesses} supplies the required factors as words in the
input subgroup generators. The centralizer and simultaneous conjugacy
algorithms are uniform in the strand number and input lists; all other
steps are explicit word operations. This proves uniformity in $n$ and
in all the supplied words.
\end{proof}

The centralizers used in the proof must be computed in the larger braid
group $B_{n+1}$;
embedding generating sets for centralizers computed in
the original braid group $B_n$ does not give the required generating sets
in the larger group.

For fixed $n$ and fixed parabolic subgroups $A,B\le B_n$, the generators
$c_1,\ldots,c_r$ and $d_1,\ldots,d_s$ of their centralizers in $B_{n+1}$
can be computed in advance. Write $|w|$ for the length of the chosen
word representative of $w$. The total length of the two unreduced input
lists in~\eqref{eq:simultaneous} is exactly
\[
2\sum_{i=1}^r|c_i|+2\sum_{j=1}^s|d_j|
+2s\bigl(|g|+|h|\bigr).
\]
Indeed, the $i$th equation in the first family contributes $2|c_i|$,
and the $j$th in the second contributes
$2|h|+|d_j|+2|g|+|d_j|$. Free cancellation can only shorten the words.
Thus, after subgroup-dependent preprocessing, the instance size and
the time needed to assemble it are $O(1+|g|+|h|)$.

\begin{corollary}\label{cor:commuting}
Given parabolic subgroups $A,B\le B_n$, braids $g,h\in B_n$, and finite
sets $U,V\subseteq B_n$, one can decide whether there are elements
$a\in A\cap C_{B_n}(U)$ and $b\in B\cap C_{B_n}(V)$
such that $h=agb$.
The algorithm is uniform in these inputs and, on positive instances,
returns words in the supplied generators of $A$ and $B$.
\end{corollary}

\begin{proof}
For braid indices $n\ge2$, add to~\eqref{eq:simultaneous} the equations
\[
\begin{aligned}
x\inv u x&=u &&(u\in U),\\
x\inv(hvh\inv)x&=gvg\inv &&(v\in V),
\end{aligned}
\]
and solve the resulting system in $B_{n+1}$.
Set $y\coloneq g\inv x\inv h$. The added equations say that $x$
centralizes $U$ and $y$ centralizes $V$. The central factors in the
decompositions $x=z_Aa$ and $y=z_Bb$ do not affect commutation.
By Proposition~\ref{prop:ambient-reduction}, there is a factorization
$h=agb$ with the required commutation properties. Conversely, such a
factorization gives a solution $x=a$.
Recover witnessing factors by adding the commutation tests to the
search of Remark~\ref{rem:witnesses}. The assertion is trivial for braid
index $n=1$.
\end{proof}

\section{Generalizations}\label{sec:discussion}

%\subsection{Other spherical Artin groups}

For a Coxeter type $X$, let $\mathcal A(X)$ denote the corresponding
Artin group. The type is \emph{spherical} if its Coxeter group is finite;
in particular, $B_n=\mathcal A(A_{n-1})$ for $n\ge2$~\cite[p.~672]{Ch92}.
Standard parabolic subgroups are generated by subsets of the standard
generators. Parabolic subgroups are their conjugates.
Gersten and Short~\cite[proof of Proposition~4.3]{GS91} proved the rationality of centralizers
of finite sets in biautomatic groups.

\begin{lemma}\label{lem:biautomatic}
There are algorithms that, given a finite presentation promised to
define a biautomatic group $G$, solve its word and simultaneous conjugacy
problems and compute finite generating sets for centralizers of finite
sets of elements. The elements are supplied as words in the presentation
generators.
\end{lemma}

\begin{proof}
Bridson and Reeves~\cite[Proposition~1.2]{BR11} uniformly construct a
biautomatic structure from the presentation. Let $L$ be its regular
language of representatives. For a word $w$, write $\overline w$ for
the group element it represents. The structure solves the word
problem~\cite[proof of Lemma~1.3]{BR11}.
Using the constructions of Gersten and Short~\cite[pp.~135--136 and proof of Lemma~4.2]{GS91},
compose the left and right multiplier automata along fixed words
$u,v$ to obtain synchronous regular relations (with the usual end padding)
\[
\begin{aligned}
R_u&\coloneq\{(p,q)\in L^2:\overline u\,\overline p=\overline q\},\\
S_v&\coloneq\{(p,q)\in L^2:\overline p\,\overline v=\overline q\}.
\end{aligned}
\]
Intersect these relations and project to the first coordinate, erasing
padding symbols. Since every element of $G$ has a representative in $L$,
the resulting regular language is
\[
L(u,v)\coloneq\{p\in L:\overline p^{-1}\overline u\,\overline p=\overline v\}.
\]
Decide simultaneous conjugacy by testing whether an accepting state is
reachable in an automaton for $L\cap\bigcap_i L(u_i,v_i)$.
For an empty system, this language is $L$.

For a finite set $U$, the intersection of the languages $L(u,u)$ over
$u\in U$ represents exactly $C_G(U)$; for $U=\varnothing$, use $L$.
Take a deterministic finite automaton accepting this language.
Remove states that do not lie on an accepting path.
For each remaining state $q$, choose paths from the initial state to
$q$ and from $q$ to an accepting state, labeled $p_q$ and $t_q$,
respectively. Use the empty word for $p_q$ at the initial state.
For every edge $q\xrightarrow{a}r$, the words $p_qat_r$ and $p_rt_r$
are accepted, so $p_qap_r^{-1}$ represents an element of $C_G(U)$.
Include $p_f$ for every accepting state $f$.
These words generate $C_G(U)$: along an accepting path ending at $f$,
the edge words telescope to its label followed by $p_f^{-1}$.
All constructions are effective in the supplied presentation and words.
\end{proof}

\begin{corollary}\label{cor:spherical}
The double coset decision and search problems for parabolic subgroups
are uniformly solvable in the Artin groups of types
$B_r$ ($r\ge2$), $D_r$ ($r\ge4$), $E_6$, $E_7$, $H_3$, and $I_2(5)$.
The Coxeter type and rank, subgroup generating lists, and two group
words are the input. On positive instances, factors are returned as
words in the supplied subgroup generators.
\end{corollary}

\begin{proof}
Use the following standard parabolic embeddings between spherical
Artin groups, obtained by restricting the Coxeter diagrams listed by
Bourbaki~\cite[Chapter~VI, \S4, no.~1, Theorems~1 and~2]{Bo68}:
\[
\begin{aligned}
\mathcal A(B_r)&\hookrightarrow\mathcal A(B_{r+1}),\\
\mathcal A(D_r)&\hookrightarrow\mathcal A(D_N),\\
\mathcal A(E_6)&\hookrightarrow\mathcal A(E_7),\\
\mathcal A(E_7)&\hookrightarrow\mathcal A(E_8),\\
\mathcal A(H_3)&\hookrightarrow\mathcal A(H_4),\\
\mathcal A(I_2(5))&\hookrightarrow\mathcal A(H_3).
\end{aligned}
\]
For type $D$, choose an even integer $N>r$: take $N\coloneq r+1$ when
$r$ is odd and $N\coloneq r+2$ when $r$ is even.
These inclusions are injective~\cite[Proposition~2.5]{DGM24} and
effective: each standard generator maps to the corresponding ambient
generator. Each source group $H$ is a proper standard parabolic
subgroup of its ambient group $G$. Every parabolic subgroup of $H$
remains parabolic in $G$, since it is a conjugate of a subgroup generated
by a subset of the same standard generators.

By the double centralizer theorem of Ajbal and
Godelle~\cite[Theorem~0.1(ii) and Proposition~2.1]{AG19},
every proper standard parabolic subgroup $Q$ of these ambient groups
satisfies
\[
C_G(C_G(Q))=Z(G)\times Q,
\]
where the product is internal. The exceptional ambient types
$D_{2k+1}$ and $E_6$ do not occur here. Taking $Q=H$ gives
$Z(G)\cap H=\{1\}$, and conjugation gives the required double centralizer
identities for the input parabolic subgroups. For the trivial subgroup,
the double centralizer is $Z(G)$.

Spherical Artin groups are biautomatic~\cite[Theorem~3.4]{Ch92}. Their standard
finite presentations are computable from their Coxeter diagrams~\cite[Introduction,
presentation~(1), and \S1.2.1]{AG19}.
Apply Lemma~\ref{lem:biautomatic} to obtain the word, simultaneous
conjugacy, and centralizer algorithms uniformly.
Proposition~\ref{prop:ambient-reduction} and Remark~\ref{rem:witnesses}
now apply.
\end{proof}

%\subsection{Finite direct products}

\begin{corollary}\label{cor:products}
Let $G\coloneq G_1\times\cdots\times G_k$, where each factor is a braid
group or one of the Artin groups in Corollary~\ref{cor:spherical}.
Equip $G$ with the union of the standard generating sets of its factors,
supported in
their respective coordinates, and define parabolic subgroups using
this generating set. The double coset decision and search problems for
parabolic subgroups of $G$ are uniformly solvable, with the factor types
and ranks included in the input.
\end{corollary}

\begin{proof}
By the choice of generators, standard parabolic subgroups are products
of standard parabolic subgroups of the factors, and conjugation acts
coordinatewise.
Thus, writing $A_i,B_i$ for the projections of $A,B$, we have
\[
\begin{aligned}
A&=\prod_{i=1}^k A_i,\\
B&=\prod_{i=1}^k B_i.
\end{aligned}
\]
Moreover, $h\in AgB$ if and only if $h_i\in A_i g_i B_i$ for all $i$.
Project the supplied generating lists and apply Main Theorem~\ref{thm:main}
and Corollary~\ref{cor:spherical} to the factors.
On positive instances, assemble the factors into tuples $a\in A$ and
$b\in B$ with $h=agb$. To express these tuples in the original supplied
subgroup generators, enumerate words and test equality coordinatewise.
The searches terminate because $a\in A$ and $b\in B$.
\end{proof}

\subsection*{Acknowledgments}
Much of this research was conducted while the first author was a
postdoctoral researcher at Bar-Ilan University, hosted by the second and
third authors. He thanks them for their hospitality during this period.
The first author acknowledges financial support from the Minerva
Foundation of Germany.
The second author acknowledges financial support from the Oswald Veblen
Fund.

\subsection*{Use of AI}
The main result was obtained in the pre-AI era.
We used ChatGPT for literature searches and proof optimization.
All proofs were formalized in Lean~4~\cite{DCPLean}.

\end{document}